\documentclass[12pt,letterpaper]{article}
\usepackage[margin=1in]{geometry}
\usepackage[T1]{fontenc}
\usepackage{amsmath,amssymb,amsthm,mathtools}
\usepackage{newtxtext,newtxmath}
\usepackage{microtype}
\usepackage{booktabs,array}
\usepackage{needspace}
\usepackage{placeins}
\usepackage{tikz}
\usetikzlibrary{calc,arrows.meta}
\usepackage[hidelinks]{hyperref}
\hypersetup{pdftitle={When Polynomial Multiplication Counts Geometry:
Schubert Polynomials and Monk's Rule},pdfauthor={Ryan Bansal}}
\newtheorem{theorem}{Theorem}[section]
\newcommand{\C}{\mathbb C}

\newcommand{\PP}{\mathbb P}
\newcommand{\Fl}{\operatorname{Fl}}
\newcommand{\Sch}{\mathfrak S}
\newcommand{\CP}{C_{\mathrm{pt}}}
\newcommand{\CL}{C_{\mathrm{line}}}
\newcommand{\PM}{\mathcal P_M}
\newcommand{\Lq}{\mathcal L_q}

\title{When Polynomial Multiplication Counts Geometry:\\
Schubert Polynomials and Monk's Rule}
\author{Ryan Bansal\\Lynbrook High School\\San Jose, CA}
\date{}

\begin{document}
\maketitle
\begin{abstract}
\normalsize
We give a geometric introduction to Schubert polynomial multiplication
through point--line flags in the complex projective plane. Moving a
prescribed point onto a prescribed line gives a family with local
equation \(au=\lambda\), whose limiting curve has two components.
The incidence equations and their local multiplicities give a
two-term divisor product, using the standard intersection product on
classes of algebraic cycles. Permutation labels and Schubert polynomials
express the same decomposition as an instance of Monk's rule. For three
incidence conditions, each allowing either a point on a prescribed line
or a line through a prescribed point, the product has degree six.
Under explicit general-position assumptions, direct construction and
local equations give six distinct solutions of multiplicity one.
\end{abstract}

\section{A point and a line}\label{sec:opening}

Two distinct points determine a line. In the projective plane, two
distinct lines also meet in a point. Enumerative geometry asks how many
configurations satisfy a given collection of conditions. To combine
those conditions systematically, we need an arithmetic in which
multiplication means imposing them at once.

This was the principle behind Schubert's nineteenth-century calculus:
geometric conditions became symbols that could be multiplied. Moving
the data into special positions often made the answers easier to see.
Justifying this method was part
of Hilbert's fifteenth problem~\cite[Problem~15]{Hilbert}; modern
intersection theory provides its foundations~\cite[\S2]{BGP}.

Point--line flags give a small model of this calculus~\cite[\S3.4]{Smirnov}.
A point \(p\) together with a line \(L\) through it is called a \emph{flag}.
Fix a line \(M\) and a point \(q_0\) on \(M\), and require
\(p\in M\) and \(q_0\in L\). If \(p\ne q_0\), the line through these
two points must be \(M\), so \(L=M\). If \(p=q_0\), any line through
that point is allowed. These two kinds of flag will become the summands in
\begin{equation}\label{eq:preview}
 \underbrace{x_1}_{p\in M}\,
 \underbrace{(x_1+x_2)}_{q_0\in L}
 =\underbrace{x_1^2}_{p=q_0}
  +\underbrace{x_1x_2}_{L=M}.
\end{equation}
The two factors record the conditions imposed together, while the
summands record the resulting kinds of flag. Here the variables
represent intersection classes, not coordinates or equations defining
the conditions. The geometry explains not only which terms appear,
but also why each has coefficient one.

The two-term decomposition suggests a way to tackle a finite counting
problem. Given
three lines \(M_i\) and three points \(q_i\), require, for each
\(i=1,2,3\),
\[
 p\in M_i\quad\text{or}\quad q_i\in L.
\]
For suitable choices of the prescribed points and lines, these
conditions have exactly six solutions. Multiplying their classes
predicts this number, but a weighted count need not count distinct
flags. The incidence equations determine when the solutions stay
separate and each contributes one.

\section{Points, lines, and flags}\label{sec:flags}

\subsection*{The complex projective plane}

Parallel lines do not meet in the affine plane. Projective geometry
includes points representing their directions, so that any two distinct
lines meet. Algebraically, a projective point is a one-dimensional subspace of
\(\C^3\). A nonzero vector \((x,y,z)\) specifies the point
\([x:y:z]\); multiplying all three entries by the same nonzero scalar
does not change it. The set of these points is the projective plane
\(\PP^2\). More generally, \(\PP(V)\) denotes the one-dimensional
subspaces of a vector space \(V\).
Where \(y\ne0\), rescaling gives \([x/y:1:z/y]\), described by the
ordinary pair \((x/y,z/y)\). The points \([x:0:z]\) supply the
missing directions. Similarly, the projective line \(\PP^1\) consists
of nonzero pairs \([s:t]\), up to nonzero rescaling: one coordinate \(s/t\)
when \(t\ne0\), together with the point \([1:0]\).

Complex numbers are needed for a similar reason: real solutions can
disappear as the data move. Consider
\[
 t^2=\varepsilon.
\]
For every nonzero complex \(\varepsilon\), there are two distinct
complex roots. Over \(\mathbb R\), they disappear when a real
\(\varepsilon\) changes from positive to negative. At \(\varepsilon=0\),
the two complex roots coincide: the equation is \(t^2=0\), with the
root zero counted twice. Counting complex solutions with weights,
called \emph{multiplicities}, keeps the total equal to two.

\subsection*{Line coordinates and flags}

A projective line is \(\PP(W)\) for a two-dimensional subspace
\(W\subset\C^3\). Such a subspace is the kernel of a nonzero linear
form, so a line has an equation \(ax+by+cz=0\).
Rescaling \((a,b,c)\) by a common nonzero scalar changes neither the
kernel nor the line. Conversely, linear forms with the same kernel
are scalar multiples. Consequently a line itself has coordinates
\([a:b:c]\), and the space of all lines is the \emph{dual projective plane}
\((\PP^2)^*\). Incidence becomes
\begin{equation}\label{eq:incidence}
 p=[x:y:z],\qquad L=[a:b:c],\qquad
 p\in L\ \Longleftrightarrow\ ax+by+cz=0.
\end{equation}
Point coordinates and line coordinates may be rescaled independently.
Two distinct projective lines meet because their two-dimensional
kernels in \(\C^3\) intersect in a one-dimensional subspace.

An incident pair is a \emph{complete flag} in \(\C^3\):
\[
 0\subset F_1\subset F_2\subset\C^3,\qquad
 \dim F_1=1,\quad \dim F_2=2,\qquad
 p=\PP(F_1),\quad L=\PP(F_2).
\]
Here the nesting means precisely that \(p\in L\).
The \emph{flag variety} is the space of all these pairs:
\[
 \Fl_3=\{(p,L)\in\PP^2\times(\PP^2)^*:p\in L\}.
\]
A point of \(\Fl_3\) is therefore one entire flag, not just its point
in the original plane. Its dimension is three: two parameters choose
\(p\), and one more chooses a line through \(p\). Throughout, dimensions
count complex parameters.

\subsection*{Local coordinates}

Where the point's second coordinate and the line's third coefficient
are nonzero, scale each independently to obtain
\[
 p=[u:1:z],\qquad L=[a:b:1].
\]
Equation~\eqref{eq:incidence} now gives
\begin{equation}\label{eq:chart}
 b=-au-z.
\end{equation}
Conversely, this formula gives a flag in the patch for every triple
\((u,a,z)\in\C^3\), so these are three independent local coordinates.
Although the patch does not contain every flag, analogous patches
together cover \(\Fl_3\). One can always choose a nonzero
point coordinate and a nonzero line coefficient in different positions;
otherwise incidence would have just one nonzero term and could not
vanish. Normalize those two entries and use incidence to solve for the
line coefficient in the position of the normalized point coordinate.
Each patch has three free coordinates. In particular,
\(\Fl_3\) is \emph{smooth}: locally it has ordinary complex coordinates,
with no singularity. It is also projective, since its incidence equation
cuts out a closed subset of a product of projective planes.

For a prescribed line \(M\) and prescribed point \(q\), write
\[
 \PM=\{(p,L):p\in M\},\qquad
 \Lq=\{(p,L):q\in L\}.
\]
Such a solution set is called a \emph{locus}. In \(\PM\), one parameter
chooses \(p\) on \(M\) and another chooses \(L\) through \(p\).
For \(\Lq\), choose \(L\) through \(q\), then \(p\) on \(L\).
Both loci have dimension two, so their \emph{codimension} in \(\Fl_3\)
is \(3-2=1\). A codimension-one algebraic piece is called a
\emph{divisor}. Imposing these two divisors together gives a curve of flags.

\section{A moving intersection}\label{sec:degeneration}

First place \(q\) off \(M\). A flag satisfying both conditions is obtained
by choosing any \(p\in M\) and joining it to \(q\). The choices form
a projective line \(\PP^1\), so their locus is a curve of flags.
Now move \(q\) onto \(M\). What happens to this curve?

Choose
\[
 M=\{z=0\},\qquad q_0=[0:1:0],\qquad
 q_\lambda=[0:1:\lambda].
\]
As \(\lambda\) varies, the prescribed point moves along \(x=0\).
This line describes the point's path; it imposes no further condition
on the flag.
In the chart~\eqref{eq:chart},
\[
 p\in M\iff z=0,\qquad
 q_\lambda\in L\iff b+\lambda=0\iff au+z=\lambda.
\]
On the surface \(z=0\), the remaining coordinates are \((u,a)\).
The intersection \(D_\lambda=\PM\cap\mathcal L_{q_\lambda}\)
therefore has equation
\begin{equation}\label{eq:family}
 au=\lambda.
\end{equation}

At \(\lambda=0\), the equation factors as \(au=0\).
The factor \(u=0\) fixes the point \(p=q_0\); the factor \(a=0\)
fixes the line \(L=M\). These are the two curve components
\begin{equation}\label{eq:components}
 \CP=\{(p,L):p=q_0\},\qquad
 \CL=\{(p,L):L=M\},\qquad
 D_0=\CP\cup\CL.
\end{equation}
A \emph{component} is an algebraic piece that cannot be split further
in this way. The incidence argument identifies the complete locus,
including flags outside the chart: if \(p\ne q_0\), then \(L=M\);
if \(p=q_0\), every line through it is allowed. Each component is a
projective line, and they meet at the flag \((q_0,M)\).

\begin{figure}[htbp]
\centering
\begin{tikzpicture}[>=Latex,every node/.style={font=\small},scale=0.82]
 \begin{scope}
  \coordinate (p) at (-2,0);
  \coordinate (q0) at (1,0);
  \coordinate (q) at (1,1.35);
  \draw[thick] (-2.7,0)--(2,0) node[right] {\(M\)};
  \draw[thick] ($(p)!-0.15!(q)$)--($(p)!1.15!(q)$)
    node[above right] {\(L_\lambda\)};
  \fill (p) circle (2pt) node[below] {\(p\)};
  \fill (q0) circle (2pt) node[below] {\(q_0\)};
  \fill (q) circle (2pt) node[above left] {\(q_\lambda\)};
  \draw[->,very thick,blue!70!black]
    (1,1.12)--(1,0.2);
 \end{scope}
 \begin{scope}[xshift=6.4cm]
  \coordinate (q0) at (0,0);
  \coordinate (q) at (0,1.35);
  \coordinate (p) at (1.7,0);
  \draw[thick] (-1.2,0)--(2.6,0) node[right] {\(M\)};
  \draw[thick] (-0.45,1.707)--(2.4,-0.556)
    node[right] {\(L_\lambda\)};
  \draw[densely dashed] (-0.95,0.754)--(1.4,-1.112)
    node[right] {\(L_0\)};
  \fill (q0) circle (2pt) node[below left] {\(q_0\)};
  \fill (q) circle (2pt) node[above left] {\(q_\lambda\)};
  \fill (p) circle (2pt) node[above right] {\(p_\lambda\)};
  \draw[->,very thick,blue!70!black] (0,1.1)--(0,0.2);
  \draw[->,very thick,blue!70!black] (1.42,0)--(0.23,0);
 \end{scope}
\end{tikzpicture}
\caption{Two limiting motions}\label{fig:moving}
\end{figure}
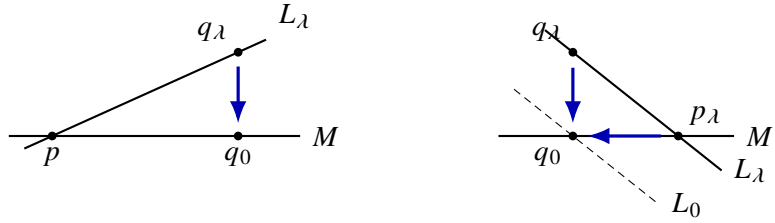

Both components occur as limits. Holding \(u=u_0\ne0\) and
taking \(a=\lambda/u_0\) approaches \(a=0\), the fixed-line component.
Holding \(a=a_0\ne0\) and taking \(u=\lambda/a_0\) approaches \(u=0\),
the fixed-point component. In the second path the flag's point moves
toward \(q_0\) as well as the prescribed point.
Figure~\ref{fig:moving} shows these two motions in the original plane.

\begin{figure}[htbp]
\centering
\begin{tikzpicture}[>=Latex,every node/.style={font=\small},scale=0.78]
 \draw[->,thick] (-2.5,0)--(2.6,0) node[right] {\(a\)};
 \draw[->,thick] (0,-2.5)--(0,2.6) node[above] {\(u\)};
 \node[below] at (1.8,0) {\(\CP\)};
 \node[left] at (0,1.8) {\(\CL\)};
 \draw[thick,blue!70!black,domain=0.25:2.2,samples=100]
   plot (\x,{0.5/\x});
 \draw[thick,blue!70!black,domain=-2.2:-0.25,samples=100]
   plot (\x,{0.5/\x});
 \node[blue!70!black] at (1.55,1.2) {\(D_\lambda\)};
 \fill (0,0) circle (2pt);
\end{tikzpicture}
\caption{Two limiting components}\label{fig:limit}
\end{figure}
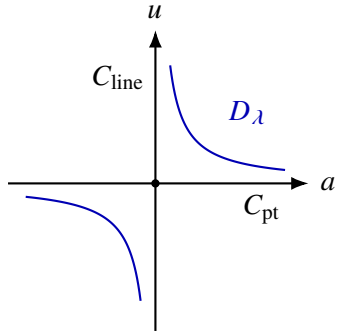

Figure~\ref{fig:limit}, in contrast, is in flag parameter space:
every point stands for an entire flag. For positive real \(\lambda\),
the drawing has two hyperbola branches; they are parts of one complex
curve, not two moving components. Globally, for \(\lambda\ne0\),
choosing \(p\in M\) still gives a unique flag, so \(D_\lambda\)
is one projective-line curve.

Passing to such a limiting configuration is called a
\emph{degeneration}. The limit has two components, but the general
curve has only one. An arithmetic of intersections must account for
this change without treating a curve as equal to a different set of flags.

\FloatBarrier
\section{Classes and intersection products}\label{sec:classes}

Different lines in \(\PP^2\) occupy different positions, but each meets
a further general line once. To keep this intersection behavior while
forgetting a particular position, we work with \emph{classes}.

\subsection*{Weighted collections and their classes}

An \emph{algebraic} locus is given by polynomial equations.
It is \emph{irreducible} if it is not the union of two smaller closed
algebraic loci. A \emph{cycle} is a finite formal integer sum of
irreducible loci of the same dimension, such as
\(\CP+\CL\) or \(2\CP\). Here \emph{formal} means that addition records
the pieces with coefficients; it is not ordinary set union.
Positive coefficients can record multiplicities; allowing negative
coefficients makes subtraction possible.

The simplest example is replacing one projective line by another.
The zero lines of any two nonzero linear forms on \(\PP^2\) will have
the same class, denoted \(H\).
In the flag space, a \emph{subvariety} is an algebraic locus.
To compare curves where equations vanish, counted with multiplicity,
we restrict the equations to the same irreducible surface.
An equation is \emph{homogeneous separately} if every term has
the same total degree in the point coordinates, and likewise in the
line coordinates. For instance, both \(b+\lambda c\) and \(b\) have
point degree zero and line degree one, before any chart normalization.
Their zero curves on \(\PM\) are the moving curve and its limit.

On an irreducible subvariety, two homogeneous equations having the same
degree separately in the point coordinates and in the line coordinates
define equivalent weighted zero loci, provided neither equation vanishes
identically on that subvariety. Two cycles have the \emph{same class}
when their difference is a finite integer sum of these zero-locus
differences. Every comparison needed in this paper uses linear
equations.\footnote{For \(k\)-cycles, the supporting subvariety \(W\)
is closed, irreducible, and of dimension \(k+1\). Each equation is
homogeneous in each coordinate block; the pair has the same two
degrees, but these may vary between pairs. Zeros are taken after
restriction to \(W\), with their divisor multiplicities. The ratio is
a nonzero rational function on \(W\), whose principal divisor is the
zero-locus difference. Conversely, on a nonempty coordinate chart any
rational function is a polynomial ratio; homogenizing numerator and
denominator to matching separate degrees gives such a pair. Thus this
is rational equivalence, not a restriction to degree-one
relations~\cite[Tag~0AZI]{Stacks}. On singular \(W\), the weights are the
standard divisor multiplicities.}

For equations \(F=0\) and \(G=0\) of matching degrees, the ratio
\(F/G\), defined where \(G\ne0\), is unchanged by rescaling either
coordinate block. Thus the comparison does not depend on chosen
coordinate vectors. We may replace either weighted zero locus by the
other inside a cycle expression without changing its class.
Square brackets, such as \([C]\), denote a cycle's class; \(C\) is
called a \emph{representative} of that class.
Multiplicity is essential: \(t^2=0\) supplies twice the zero of \(t=0\),
despite having the same underlying set.

On \(\Fl_3\), linear forms in the point coordinates compare the conditions that \(p\)
lie on different prescribed lines. Linear forms in the line coordinates
compare the conditions that \(L\) contain different prescribed points.
None of these forms vanishes on every flag. Consequently
\begin{equation}\label{eq:divisor-classes}
 [\mathcal P_M]=A\quad\text{for every line }M,\qquad
 [\mathcal L_q]=B\quad\text{for every point }q.
\end{equation}

\subsection*{Intersection products}

On a smooth projective variety these classes have an associative,
commutative, distributive multiplication, the \emph{intersection product}.
Changing a representative within its class does not change the product.
The product imposes conditions together. For suitably intersecting
representatives, it is their intersection, with its local
multiplicities.\footnote{The classes form the Chow ring
\(\mathrm{CH}^*(\Fl_3)\). Proper representatives compute the product with
intersection multiplicities~\cite[Tag~0B0G]{Stacks}. In our divisor
calculations, a simple local cut gives multiplicity one; the divisor
multiplicity is described in~\cite[Tag~0B05]{Stacks}.}
Products add codimensions. The class \(1=[\Fl_3]\) imposes no condition;
the class \(0\) can be represented by an empty locus.

An intersection is \emph{proper} if every component has the expected
codimension. In this three-dimensional space, two divisor conditions
should give curves, and three should give isolated flags. Properness
alone does not make the weights one. Move the chosen flag to the origin
in local coordinates and take the degree-one terms of the equations.
If these \emph{linear parts} are independent, the cut is \emph{simple}
and has multiplicity one.
For example, \(t=0\) has nonzero linear part \(t\), whereas
\(t^2=0\) has zero linear part. Independence of several linear parts
is checked by row reduction.
For curve components, simplicity at ordinary points of each component
is enough to determine their coefficients; exceptional crossings
do not determine the weight of the whole curve.

Changing representatives is particularly important for squares.
The product \(H^2\) uses two distinct lines, whose intersection is
one point, not the whole line obtained by intersecting a line with
itself as a set. Three lines with no common point give \(H^3=0\).
Similarly, \(A^2\) uses two different prescribed lines by
\eqref{eq:divisor-classes}. A \emph{general} choice means positions
avoiding unwanted coincidences; each calculation below specifies
what must be avoided.

\subsection*{The class of the moving curve}

On the fixed surface \(\PM\), a flag has \(p=[u:v:0]\),
\(L=[a:b:c]\), with \(au+bv=0\). This surface is
irreducible.\footnote{Here we use the standard facts that
\(\C\times\PP^1\) is irreducible and that overlapping irreducible
open patches have irreducible union. Where \(u\ne0\), normalize
\(p=[1:v:0]\); then \(a=-bv\), leaving free choices
\(v\in\C\) and \([b:c]\in\PP^1\). This patch is
\(\C\times\PP^1\). The analogous patch where \(v\ne0\) overlaps
it, and the two cover \(\PM\).}
The equations
\[
 b+\lambda c=0,\qquad b=0
\]
have the same degree in the line coordinates and degree zero in the
point coordinates. Neither vanishes identically on \(\PM\).
Both equations lie in the linear family \(sb+tc\), for
\([s:t]\in\PP^1\).
The restrictions of \(b,c\) are independent, so no member vanishes
identically. Their weighted zero curves therefore have the same class.

For \(\lambda\ne0\), the first zero curve is \(D_\lambda\).
In chart~\eqref{eq:chart}, its equation is \(au=\lambda\).
At a point \((u_0,a_0)\) on it, \(a_0u_0=\lambda\), so at least one
coefficient of its linear part
\[
 a_0(u-u_0)+u_0(a-a_0)
\]
is nonzero. Hence the curve has multiplicity one.
The equation \(b=0\) becomes \(au=0\).
Along \(u=0\) away from the crossing, \(a\ne0\), so the factor \(u\)
occurs once. Along \(a=0\), the same argument applies to \(a\).
This chart contains ordinary parts of every curve in question; the
global incidence argument already identified all components.
It follows that
\[
 [D_\lambda]=[\CP]+[\CL].
\]
The original two divisors intersect properly. Write \(u=u_0+h\),
\(a=a_0+k\), \(z=l\) at a flag on \(D_\lambda\).
The equations \(z=0\), \(au+z=\lambda\) have linear parts
\(l\) and \(a_0h+u_0k+l\); the remaining term is \(hk\).
Their coefficient rows \((0,0,1)\) and \((a_0,u_0,1)\) are
independent because \(a_0u_0=\lambda\ne0\).
Thus their intersection computes \(AB\), giving
\begin{equation}\label{eq:geometric-product}
 \boxed{AB=[\CP]+[\CL].}
\end{equation}
The special intersection confirms the same sum: \(z=0\), \(au=0\)
gives the two components once each. They fail to cut simply at
their crossing, but coefficients of curve components are determined
away from that crossing. Their common flag is not a third curve
component, and no point is subtracted from a sum of curves.

\subsection*{Squares and top products}

Choose another line \(N\ne M\) through \(q_0\). The conditions
\(p\in M\) and \(p\in N\) fix \(p=q_0\) and leave \(L\) free.
Two independent line equations cut the point simply, so
\(A^2=[\CP]\). Dually, two distinct prescribed points on \(M\)
force \(L=M\) and leave \(p\) free, giving \(B^2=[\CL]\). Therefore
\begin{equation}\label{eq:squares}
 A^2=[\CP],\qquad B^2=[\CL],\qquad AB=A^2+B^2.
\end{equation}
A third line not through \(q_0\) leaves no point, so \(A^3=0\).
Three noncollinear prescribed points leave no line, so \(B^3=0\).

Take \(q'\in M\setminus\{q_0\}\). Imposing \(q'\in L\) on \(\CP\)
leaves just \((q_0,M)\). Imposing \(p\in N\) on \(\CL\) leaves the
same flag. Each additional linear equation cuts its projective-line
family at one simple point. Write
\begin{equation}\label{eq:top-products}
 T=[(q_0,M)],\qquad A^2B=AB^2=T,\qquad A^3=B^3=0.
\end{equation}
The class \(T\) is represented by one isolated flag, whereas
\(A^2\) and \(B^2\) are represented by curves of flags.

For a zero-dimensional cycle over \(\C\), its \emph{degree} is the
sum of its coefficients, including multiplicities. On a projective
variety this sum is unchanged by class equivalence~\cite[Tag~0AZ0]{Stacks}.
In particular, \(\deg T=1\). A zero-dimensional intersection
representing \(dT\) has total weight \(d\), and consists of \(d\)
distinct flags when every local multiplicity is one.

These products also determine the coefficients of a curve expression.
If \(\gamma=aA^2+bB^2\), then~\eqref{eq:top-products} gives
\begin{equation}\label{eq:coefficient-tests}
 \deg(\gamma B)=a,\qquad \deg(\gamma A)=b.
\end{equation}
Geometrically, requiring the flag's line to contain a further general
point selects one flag on \(\CP\) and none on \(\CL\); requiring
its point to lie on a further general line does the reverse.
The coefficients describe intersection behavior, not the number of
components: the moving curve \(D_\lambda\) is one irreducible curve
with class \(A^2+B^2\).

\section{Six positions and their polynomials}\label{sec:positions}

The two curves describe special positions relative to the reference
flag \(q_0\in M\). To find all the possible positions, test
\[
 p=q_0,\qquad p\in M,\qquad q_0\in L,\qquad L=M.
\]
These tests are not independent. Either equality forces both middle
incidences; if both incidences hold and \(p\ne q_0\), then \(L=M\).
The six exact positions are: neither incidence; only \(p\in M\);
only \(q_0\in L\); \(p=q_0\) with \(L\ne M\);
\(L=M\) with \(p\ne q_0\); and \((p,L)=(q_0,M)\).
Every flag belongs to exactly one of them.

An exact position need not include its limits. Its \emph{closure} adds
those limiting configurations. For example, flags with only \(p\in M\)
can approach flags with \(p=q_0\) or \(L=M\); their closure is
\(\PM\). The exact positions are called \emph{Schubert cells},
and their closures are \emph{Schubert varieties}. Unlike the cells,
the closures can overlap.

\subsection*{Permutation labels}

Permutation labels encode these positions so that multiplication can
be computed by swapping entries, without solving new incidence equations.
A permutation is an ordering of \(1,2,3\); the set of these orderings
is denoted \(S_3\). In \emph{one-line notation}, \(213\) means that
the first value is 2, the second 1, and the third 3.
An \emph{inversion} is an out-of-order pair: positions \(i<j\)
with \(w(i)>w(j)\). The number of inversions is the \emph{length} \(\ell(w)\).
For example, \(213\) has one inversion, while \(312\) has two.
Our convention makes \(\ell(w)\) the codimension of the position
labelled by \(w\).

The labels can be read from a basis adapted to both flags. Suppose only
\(p\in M\) holds, so \(p\ne q_0\) and \(q_0\notin L\).
Choose \(v_3\) along \(q_0\), \(v_2\) along \(p\), and
\(v_1\) in the vector plane for \(L\), outside the plane for \(M\).
These vectors form a basis. Numbering backward along the reference
flag, its subspaces and those of the moving flag are
\[
 \begin{aligned}
 E_1&=\langle v_3\rangle,&
 E_2&=\langle v_3,v_2\rangle,& E_3&=\C^3,\\
 F_1&=\langle v_2\rangle,&
 F_2&=\langle v_2,v_1\rangle,& F_3&=\C^3.
 \end{aligned}
\]
Reading the order in which directions enter the moving flag gives
\(2,1,3\). The last digit records the remaining direction needed for
the full space. Thus this position is \(213\).
This basis is adapted to the particular pair of flags, rather than
fixed once for all pairs.

\begin{figure}[htbp]
\centering
\begin{tikzpicture}[
 every node/.style={font=\small,align=center,inner sep=3pt},
 edge/.style={draw=black,shorten >=4pt,shorten <=4pt}]
 \node (top) at (0,0) {\((p,L)=(q_0,M)\)\\{\scriptsize\(321\)}};
 \node (P) at (-2.35,-1.3) {\(p=q_0\)\\{\scriptsize\(312\)}};
 \node (Q) at (2.35,-1.3) {\(L=M\)\\{\scriptsize\(231\)}};
 \node (A) at (-2.35,-2.7) {\(p\in M\)\\{\scriptsize\(213\)}};
 \node (B) at (2.35,-2.7) {\(q_0\in L\)\\{\scriptsize\(132\)}};
 \node (one) at (0,-4) {\(\Fl_3\)\\{\scriptsize\(123\)}};
 \node[anchor=east,font=\scriptsize] at (-3.65,0) {codim \(3\)};
 \node[anchor=east,font=\scriptsize] at (-3.65,-1.3) {codim \(2\)};
 \node[anchor=east,font=\scriptsize] at (-3.65,-2.7) {codim \(1\)};
 \node[anchor=east,font=\scriptsize] at (-3.65,-4) {codim \(0\)};
 \draw[edge] (one)--(A); \draw[edge] (one)--(B);
 \draw[edge] (A)--(P); \draw[edge] (A)--(Q);
 \draw[edge,preaction={draw=white,line width=3pt}] (B)--(P);
 \draw[edge] (B)--(Q);
 \draw[edge] (P)--(top); \draw[edge] (Q)--(top);
\end{tikzpicture}
\caption{Six Schubert closures}\label{fig:six}
\end{figure}
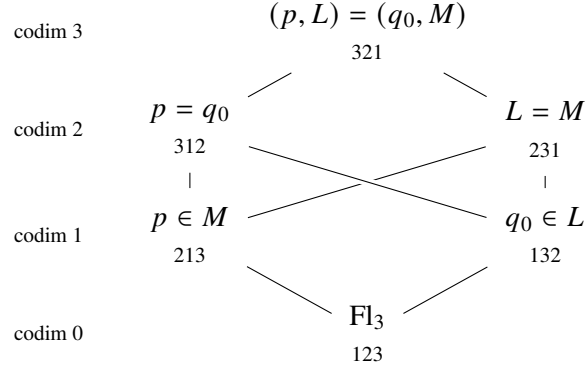

Applying the same numbering to \(p=q_0\), \(L\ne M\), direction 3
enters first, then a direction outside \(M\), numbered 1, then 2.
This gives \(312\). For \(L=M\), \(p\ne q_0\), the order is
2, then 3, then 1, giving \(231\).
The same convention gives \(123\) in generic position and \(321\)
for the reference flag. Write \(\Omega_w\) for the closure with label
\(w\). The standard classification supplies Table~\ref{tab:dictionary},
with
\(\operatorname{codim}\Omega_w=\ell(w)\)~\cite[Proposition~1.2.1
and subsequent cell description]{Brion}.\footnote{In the
dimension-indexed notation of Brion and of the geometric varieties in
Billey--Gao--Pawlowski, our \(\Omega_w=X_{w_0w}\), where \(w_0=321\).
The composition \(w_0w\) reverses values, not positions. Numbering
the adapted basis backward implements this conversion.}
Figure~\ref{fig:six} records their containments.

Each edge joins a closure to a contained closure of codimension
one higher. These one-step containments are called covers in
\emph{Bruhat order}. The diagram is unweighted: it records which
closures contain which, not the coefficient of a multiplication.

\subsection*{Polynomial representatives}

The arithmetic in \(A,B\) describes this flag space. Schubert
polynomials place it in a family of polynomial representatives indexed
by permutations, with the same construction in larger flag spaces.
Introduced by Lascoux and Sch\"utzenberger~\cite{LS}, they give a
standard polynomial \(\Sch_w\) representing \([\Omega_w]\).
The standard choices for \(A,B\) are \(x_1,x_1+x_2\), respectively.
To read these polynomials as classes, substitute
\begin{equation}\label{eq:substitution}
 x_1\longmapsto A,\qquad x_2\longmapsto B-A.
\end{equation}
Then \(x_1+x_2\) represents \(B\). The variable \(x_2\) alone
represents a difference of classes, not another incidence condition.
These are the general divisor coordinates of~\cite[Lemma~3.114]{BGP}.

In Table~\ref{tab:dictionary}, the polynomial degree matches the
codimension in every row. The first column describes the closures,
not the six disjoint exact positions.

\begin{table}[htbp]
\centering
\renewcommand{\arraystretch}{1.18}
\begin{tabular}{@{}lcccc@{}}
\toprule
Closure condition & \(w\) & Codim. & Class & \(\Sch_w\)\\
\midrule
all flags & \(123\) & \(0\) & \(1\) & \(1\)\\
\(p\in M\) & \(213\) & \(1\) & \(A\) & \(x_1\)\\
\(q_0\in L\) & \(132\) & \(1\) & \(B\) & \(x_1+x_2\)\\
\(p=q_0\) & \(312\) & \(2\) & \(A^2\) & \(x_1^2\)\\
\(L=M\) & \(231\) & \(2\) & \(B^2\) & \(x_1x_2\)\\
\((p,L)=(q_0,M)\) & \(321\) & \(3\) & \(T\) & \(x_1^2x_2\)\\
\bottomrule
\end{tabular}
\caption{The polynomial dictionary}\label{tab:dictionary}
\end{table}

The first four rows follow immediately from~\eqref{eq:substitution}.
For the remaining rows, the geometric products give
\[
 x_1x_2\longmapsto AB-A^2=B^2,\qquad
 x_1^2x_2\longmapsto A^2B-A^3=T.
\]
In particular, a polynomial representing a class is not an equation
defining its geometric locus.

\subsection*{Divided differences}

Many polynomials can represent the same class. The standard Schubert
polynomials are obtained from \(\Sch_{321}=x_1^2x_2\) by
\emph{divided differences}.
Let \(s_i f\) mean the polynomial obtained by exchanging
\(x_i,x_{i+1}\) in \(f\). Define
\[
 \partial_i f=\frac{f-s_i f}{x_i-x_{i+1}}.
\]
The numerator vanishes when \(x_i=x_{i+1}\), so the division is exact
and the result is again a polynomial.
For labels, \(ws_i\) instead means swapping the entries in positions
\(i,i+1\); thus \(321s_1=231\).
When this swap decreases length by one, the standard recursion gives
\(\Sch_{ws_i}=\partial_i\Sch_w\)~\cite[Eq.~(3.50)]{BGP}.

For example,
\begin{align}
 \Sch_{231}
 &=\frac{x_1^2x_2-x_2^2x_1}{x_1-x_2}=x_1x_2,
 \label{eq:dd-line}\\
 \Sch_{312}
 &=\frac{x_1^2x_2-x_1^2x_3}{x_2-x_3}=x_1^2.
 \label{eq:dd-point}
\end{align}
The second computation temporarily uses a third formal variable
\(x_3\) to perform the exchange. It cancels from the answer; all six
polynomials in our dictionary use only \(x_1,x_2\).
For the next step, \(312s_1=132\), and
\(\partial_1(x_1^2)=(x_1^2-x_2^2)/(x_1-x_2)=x_1+x_2\).
This explains the standard sum representing \(B\); continuing gives
the remaining rows of Table~\ref{tab:dictionary}.

\subsection*{Equality of polynomials versus equality of classes}

Our first product is an ordinary polynomial identity:
\begin{equation}\label{eq:schubert-product}
 \Sch_{213}\Sch_{132}=\Sch_{312}+\Sch_{231}.
\end{equation}
The two summands represent the two curves of the degeneration,
each once.
Some products require replacements to express their classes in the
six-row dictionary.

Write \(f\equiv g\) when the two polynomials give the same class under
\eqref{eq:substitution}; this is not ordinary polynomial equality.
The geometric facts imply
\begin{equation}\label{eq:replacements}
 x_1^3\equiv0,\qquad x_2^2\equiv-x_1^2-x_1x_2.
\end{equation}
Indeed, the first is \(A^3=0\). For the second, the polynomial
\(x_1^2+x_1x_2+x_2^2\) substitutes to
\[
 A^2+A(B-A)+(B-A)^2=A^2+B^2-AB=0.
\]
Substitution preserves addition and multiplication. Thus a valid
replacement can be used inside a larger product.
For example,
\begin{equation}\label{eq:square-reduction}
 (x_1+x_2)^2=x_1^2+2x_1x_2+x_2^2\equiv x_1x_2.
\end{equation}
After replacement, the product represents \(B^2\), the curve with \(L=M\),
once.
The coefficient two in the unreduced monomial expansion is not
a geometric multiplicity.

\FloatBarrier
\section{Monk's rule}\label{sec:monk}

Divided differences construct polynomial representatives; Monk's rule
computes multiplication of classes by a divisor. The divisor labels
\(213=s_1\) and \(132=s_2\) swap adjacent positions in \(123\).
The index of that swap determines a cut: after position 1 for
\(A=[\Omega_{213}]\), and after position 2 for \(B=[\Omega_{132}]\).
Across that cut, swap two entries and retain the result only when
the inversion count increases by exactly one. Each retained class
has coefficient one.

For \(AB\), multiply \(213\) by \(B\), placing the cut in \(21\mid3\).
The cross-cut swaps give
\[
 213\,t_{13}=312,\qquad 213\,t_{23}=231.
\]
Here \(t_{ij}\) swaps entries in positions \(i,j\).
Both resulting permutations have length two, one more than
\(\ell(213)=1\). Hence the rule gives
\[
 AB=[\Omega_{312}]+[\Omega_{231}],
\]
the two terms already obtained geometrically.

The length condition can discard a cross-cut swap that increases
length by too much: \(123\,t_{13}=321\) changes length from zero to
three, not one, and is excluded. Adding one divisor adds one to
codimension, which explains the one-step requirement.
It can also discard every swap. For \(A^2A\), use \(3\mid12\);
its two swaps give \(132\) and \(213\), both of smaller length than
\(312\). Neither survives, so \(A^3=0\).
An empty sum gives the zero class.

The square \(B^2\) gives a third check. From \(13\mid2\), the swap
of positions \(1,3\) gives \(231\), of length two; the swap of
positions \(2,3\) gives \(123\), of length zero. Only the first
is retained, so \(B^2=[\Omega_{231}]\). This agrees both with
the geometry and with~\eqref{eq:square-reduction}.

\subsection*{Complete flags and the general rule}

A complete flag in \(\C^n\) is a chain \(F_i\) of subspaces of
dimension \(i\), for \(i=0,\ldots,n\), with \(F_0=0\) and
\(F_n=\C^n\). Fix a reference flag \(E_i\).
Any two complete flags admit a common
adapted basis: each flag
is obtained by adding those basis directions in some order.
Number them backward along the reference flag, as in our example.
The moving order gives a permutation \(w\in S_n\). Although the
basis need not be unique, the label is determined by the relative
positions of the two flags. The closure of the position labelled
by \(w\), denoted \(\Omega_w(E)\), has
codimension \(\ell(w)\)~\cite[\S1.2]{Brion}.
Complete flag varieties are smooth and projective, and their Schubert
classes are independent of the reference flag~\cite[\S1.2]{Brion},~\cite[\S3.7]{BGP}.

The simple adjacent swap \(s_r=t_{r,r+1}\) labels a divisor.
Its geometric condition is
\begin{equation}\label{eq:general-divisor}
 \Omega_{s_r}(E)=\{F:F_r\cap E_{n-r}\ne\{0\}\},
 \qquad 1\le r<n.
\end{equation}
For \(n=3,r=1\), this says the flag's point lies on the reference
line; for \(r=2\), its line contains the reference point.
Thus these are precisely \(A\) and \(B\).

\begin{theorem}[Monk's rule]\label{thm:monk}
For \(w\in S_n\) and \(1\le r<n\),
\begin{equation}\label{eq:monk}
 [\Omega_w][\Omega_{s_r}]
 =\sum_{\substack{1\le i\le r<j\le n\\
                 \ell(wt_{ij})=\ell(w)+1}}
  [\Omega_{wt_{ij}}].
\end{equation}
\end{theorem}

Monk's theorem~\cite{Monk} is stated here in the codimension
convention above.\footnote{The finite class formula is
\cite[Theorem~3.92]{BGP}, translated from dimension indexing by
\(\Omega_w=X_{w_0w}\), where \(w_0\) reverses all values.
Brion's Proposition~1.4.3 gives the more general divisor formula;
the fundamental-weight case gives the same cross-cut coefficients.
The stable polynomial identity in \cite[Theorem~3.116]{BGP} is a
different statement. Equation~\eqref{eq:monk} is an equality of
classes, not necessarily an ordinary polynomial identity.}
The degeneration proves the two-term instance for \(AB\);
the general formula, whose proof we omit, extends the calculation to
larger flag spaces. Only positions \(1,\ldots,n\) are available for
swapping, so an empty sum gives zero in the finite flag variety.

\section{Six flags}\label{sec:count}

Choose lines \(M_1,M_2,M_3\) and points \(q_1,q_2,q_3\).
For a flag \(p\in L\), impose the three conditions
\begin{equation}\label{eq:or-conditions}
 p\in M_i\quad\text{or}\quad q_i\in L,\qquad i=1,2,3.
\end{equation}
The disjunction is inclusive: a flag may satisfy both branches.
The locus for condition \(i\) is
\(\mathcal P_{M_i}\cup\mathcal L_{q_i}\).
Writing the branch equations as \(f_i=0\), \(g_i=0\), the union has
equation \(f_i g_i=0\). Each branch occurs once, so its class is
\(A+B\). The point and line data need not form prescribed flags.

\subsection*{The degree of the intersection}

For the three ordered conditions, write a word with \(A\) for the
point-on-line branch and \(B\) for the line-through-point branch.
Distributing the three factors gives eight words:
\begin{equation}\label{eq:eight-words}
 \begin{aligned}
 (A+B)^3={}&AAA+(AAB+ABA+BAA)\\
           &+(ABB+BAB+BBA)+BBB\\
          ={}&0+3T+3T+0=6T.
 \end{aligned}
\end{equation}
The words retain the order of the geometric choices even though
class multiplication is commutative. Every mixed word is \(T\)
by~\eqref{eq:top-products}. The two pure words vanish.
Since \(\deg T=1\), the predicted intersection degree is six.
Special data could still make solutions coincide or leave a
positive-dimensional locus, so the degree alone does not give the
number of distinct solutions.

In polynomial language the same prediction is
\[
 (2x_1+x_2)^3\equiv6x_1^2x_2,
 \qquad 2x_1+x_2=\Sch_{213}+\Sch_{132}.
\]
Multiplying the replacements~\eqref{eq:replacements} and applying
them again gives \(x_1x_2^2\equiv-x_1^2x_2\) and \(x_2^3\equiv0\).
The cubic expansion leaves
\(12x_1^2x_2-6x_1^2x_2=6x_1^2x_2\).

\subsection*{General position}

If a prescribed-line vertex lies on a join of two prescribed points,
one flag can satisfy several branch words. To keep the constructions
separate, assume:
\begin{enumerate}
 \item The lines are distinct and nonconcurrent: they do not all
 pass through one point.
 \item The points are distinct and noncollinear: no line contains
 all three.
 \item No \(q_j\) lies on any \(M_i\).
 \item No vertex \(M_i\cap M_j\), for \(i\ne j\), lies on a line
 joining two of the prescribed points.
\end{enumerate}
Such choices exist: choose three nonconcurrent lines, then choose the
points successively, avoiding those lines and the joins of earlier
points to the three vertices. For the third point, also avoid the join
of the first two. Only finitely many lines are forbidden at each step.
These are sufficient assumptions, not a claim that every excluded
configuration has a different answer.

\subsection*{The six solutions}

The word \(AAA\) would put \(p\) on all three lines, impossible
by nonconcurrence. The word \(BBB\) would put all three points
on \(L\), impossible by noncollinearity.
The other six words have one of two forms.

For two \(A\) choices at indices \(i,j\) and a \(B\) choice at the
remaining index \(k\), the flag is
\begin{equation}\label{eq:aab-construction}
 p=M_i\cap M_j,\qquad L=\overline{p q_k}.
\end{equation}
The bar denotes the line joining two distinct points; they are
distinct because \(q_k\) is on none of the prescribed lines.
For two \(B\) choices at \(i,j\) and an \(A\) choice at \(k\), it is
\begin{equation}\label{eq:abb-construction}
 L=\overline{q_i q_j},\qquad p=L\cap M_k.
\end{equation}
The join is well defined, and cannot equal \(M_k\), since neither
\(q_i\) nor \(q_j\) lies on \(M_k\).
Table~\ref{tab:solutions} matches the six words with these constructions.

\begin{table}[htbp]
\centering
\renewcommand{\arraystretch}{1.2}
\begin{tabular}{@{}ccc@{}}
\toprule
Word & Point \(p\) & Line \(L\)\\
\midrule
\(AAB\) & \(M_1\cap M_2\) & \(\overline{p q_3}\)\\
\(ABA\) & \(M_1\cap M_3\) & \(\overline{p q_2}\)\\
\(BAA\) & \(M_2\cap M_3\) & \(\overline{p q_1}\)\\
\(ABB\) & \(L\cap M_1\) & \(\overline{q_2 q_3}\)\\
\(BAB\) & \(L\cap M_2\) & \(\overline{q_1 q_3}\)\\
\(BBA\) & \(L\cap M_3\) & \(\overline{q_1 q_2}\)\\
\bottomrule
\end{tabular}
\caption{Six constructed flags}\label{tab:solutions}
\end{table}

Each construction satisfies the three conditions.
Conversely, every solution of~\eqref{eq:or-conditions} satisfies
at least one branch at each index. Choosing such branches gives
one of the eight words. The pure words are impossible and each
mixed word forces its displayed construction. Thus the list covers
all solutions.

\subsection*{Distinctness}

In~\eqref{eq:aab-construction}, the vertex \(p\) lies on just the two
selected lines, not the third. The line through \(p,q_k\) cannot
contain \(q_i\) or \(q_j\): that would place the vertex on a join
of two prescribed points, contrary to assumption 4.
Also \(p\notin M_k\), so at index \(k\) only the \(B\) branch holds.
At each of the other two indices only the \(A\) branch holds.

In~\eqref{eq:abb-construction}, \(L\) contains just the selected two
points, not the third. Its point on \(M_k\) lies on neither other
prescribed line, since that would place a prescribed-line vertex on
\(\overline{q_i q_j}\). At index \(k\) only the \(A\) branch holds,
and at the other indices only the \(B\) branch holds.
Every constructed flag therefore determines exactly one branch word.
Different words cannot give the same flag: the six flags are distinct.

\subsection*{Multiplicity one}

At a constructed flag, the unused factor in \(f_i g_i=0\) is nonzero
by the distinctness argument. The linear part of the product is therefore
the selected factor's linear part multiplied by that nonzero value.
Multiplying rows by nonzero numbers preserves independence. Thus
independent selected linear parts prove simplicity for the original
union equations as well.

For~\eqref{eq:aab-construction}, two independent line equations fix
the point \(p\) at a vertex. Among lines through that point, the
third condition, containing \(q_k\ne p\), is one nonzero linear
condition and selects a single line. Locally two coordinates
choose the point and a third chooses the line through it; the three
equations have independent linear parts: the first two fix the two
point coordinates, and the third has a nonzero coefficient in the
remaining line coordinate. The isolated flag has
multiplicity one. Dually, in~\eqref{eq:abb-construction}, two
conditions fix the line and the third fixes its point simply.

The intersection consists of six distinct flags, each of weight one.
Its degree is therefore its ordinary number
of solutions, six, as predicted by~\eqref{eq:eight-words}.

\subsection*{Weighted paths in Monk's rule}

Starting at \(123\), apply the divisor sum \(A+B\) three times.
An allowed position swap contributes once for each of the two cuts
it crosses. For example, \(123\to213\to312\to321\) uses the swaps
\(t_{12},t_{13},t_{23}\). These cross only cut 1, both cuts, and
only cut 2, respectively. The divisor choices are therefore \(A\),
then \(A\) or \(B\), then \(B\), giving the words \(AAB\) and \(ABB\).
Table~\ref{tab:paths} records all four paths.

\begin{table}[htbp]
\centering
\renewcommand{\arraystretch}{1.2}
\begin{tabular}{@{}lcc@{}}
\toprule
Permutation path & Divisor words & Contribution\\
\midrule
\(123\to213\to312\to321\) & \(AAB,\ ABB\) & \(2\)\\
\(123\to213\to231\to321\) & \(ABA\) & \(1\)\\
\(123\to132\to312\to321\) & \(BAB\) & \(1\)\\
\(123\to132\to231\to321\) & \(BAA,\ BBA\) & \(2\)\\
\bottomrule
\end{tabular}
\caption{Weighted permutation paths}\label{tab:paths}
\end{table}

The words record the ordered divisor choices, totaling \(2+1+1+2=6\).
Each weight-two path records two choices, not a multiplicity-two flag:
each individual Monk summand still has coefficient one.
These are class paths, not flag trajectories; their intermediate
reference loci need not be the successive actual intersections.

\subsection*{A larger flag space}

On \(\Fl_4\), a flag is \(p\in L\subset\Pi\) in projective three-space.
The three divisors in~\eqref{eq:general-divisor} require, respectively,
a point on a prescribed plane, a line meeting a prescribed line,
and a plane containing a prescribed point.
Set \(D=\sum_{r=1}^3[\Omega_{s_r}]\).
Choosing \(p\), then \(L\) through \(p\), then \(\Pi\) through \(L\)
uses \(3,2,1\) parameters, so \(\dim\Fl_4=6\).
The variety \(\Omega_{4321}\) is the single reference flag. Finite Monk multiplication
gives\footnote{A penultimate checkpoint is
\(D^5=220[\Omega_{3421}]+280[\Omega_{4231}]+220[\Omega_{4312}]\).
Each of these has one final allowed cut and swap, so the last
coefficient is \(220+280+220=720\). See also
\cite[Exercise~3.98]{BGP}.}
\[
 D^6=720[\Omega_{4321}],\qquad \deg D^6=720.
\]
This is a degree calculation in \(\Fl_4\); the polynomial replacements
for \(\Fl_3\) do not carry over. A count of 720 distinct flags would
require further geometric checks.

\Needspace{7\baselineskip}
\section*{Acknowledgements}
The author thanks Liam Spilker for guidance and feedback during the
development of this article, and Simon Rubinstein-Salzedo for lectures
on mathematical writing that influenced its expository approach.
The author also thanks Dora Woodruff for introducing them to Schubert
polynomials and for mentorship in their research on the related forest polynomials.

\end{document}